\documentclass[12pt]{article}
\usepackage{amssymb,amsmath}
\usepackage{amsthm}

\makeatletter
\newfont{\footsc}{cmcsc10 at 8truept}
\newfont{\footbf}{cmbx10 at 8truept}
\newfont{\footrm}{cmr10 at 10truept}
\renewcommand{\ps@plain}{%
\renewcommand{\@oddfoot}{\footsc the electronic journal of combinatorics
  {\footbf 12(1)} (2005), \#R62\hfil\footrm\thepage}}
\makeatother
\newtheorem{theorem}{Theorem} 
  
\newtheorem{corollary}[theorem]{Corollary} 
\newtheorem{proposition}[theorem]{Proposition}

\theoremstyle{definition}
\newtheorem{definition}[theorem]{Definition} 
\newtheorem{example}[theorem]{Example}
\newtheorem{remark}[theorem]{Remark}
\newtheorem{question}[theorem]{Question}
\numberwithin{theorem}{section}
\numberwithin{equation}{section}

\newcommand{\N}{{\bf N}}
\newcommand{\Q}{{\bf Q}}

\newcommand{\Z}{{\bf Z}}
\newcommand{\g}{\gamma}
\newcommand{\eps}{\varepsilon}
\newcommand{\ra}{\rightarrow}
\newcommand{\fq}{\mathbf{F}_{q}}
\newcommand{\vq}{\mathbb{V}_q}
\newcommand{\B}{\mathbb{B}}
\newcommand{\qsp}{\textit{q}-species }
\newcommand{\gl}{\mathrm{GL}}
\newcommand{\glnq}{\gl_n (q)}
\newcommand{\End}{\mathrm{End}}
\newcommand{\Aut}{\mathrm{Aut}}
\newcommand{\Hom}{\mathrm{Hom}}
\newcommand{\z}{\mathcal{Z}}
\newcommand{\A}{\mathbb{A}}
\newcommand{\qbinom}[2]{\genfrac{(}{)}{0pt}{}{#1}{#2}_{\! q}}
\newcommand{\wt}{\widetilde}
\newcommand{\wh}{\widehat}
\newcommand{\fix}{\text{fix}\,}
\newcommand{\Fix}{\text{Fix}\,}
\newcommand{\sym}[1]{{\{ #1 \}}}
\newcommand{\be}{\begin{equation}}
\newcommand{\ee}{\end{equation}}

\title{
   \textbf{An Introduction to \textit{q}-Species}
}
\author{
    Kent E. Morrison \\
   \small  Department of Mathematics \\[-0.8ex]
    \small California Polytechnic State University \\[-0.8ex]
    \small San Luis Obispo, CA 93407 \\[-0.8ex]
     \small \texttt{kmorriso@calpoly.edu 
}
   }
    \date{\small
  Submitted: Oct. 22, 2005; Accepted:  Nov. 19, 2005; Published: Nov. 25, 2005 \\
  \small Mathematics Subject Classifications: 05A15; 05A30, 15A33, 18B99}
\begin{document}
\maketitle

\begin{abstract} The combinatorial theory of species developed by Joyal provides a foundation for enumerative combinatorics of objects constructed from finite sets. In this paper we develop an analogous theory for the enumerative combinatorics of objects constructed from vector spaces over finite fields. Examples of these objects include subspaces, flags of subspaces, direct sum decompositions, and linear maps or matrices of various types. The unifying concept is that of a ``\textit{q}-species,'' defined to be a functor from the category of finite dimensional vector spaces over a finite field with to the category of finite sets.
\end{abstract}

\section{Definitions}
The combinatorial theory of species originated in the work of Joyal \cite{Joyal81} in 1981 and has developed into a mature theory for understanding classical enumerative combinatorics and its generating functions \cite{BLL94}. More than thirty years ago Goldman and Rota \cite{GoldmanRota69,GoldmanRota70} began the systematic exploration of the ``subset-subspace'' analogy, the foremost example being the analogy between the binomial coefficients, which count subsets,  and the $q$-binomial coefficients, which count subspaces. Their work has an interesting prehistory that is outlined in a short survey by Kung \cite{Kung95}. The aim of this paper is to further the subset-subspace analogy with the development of the theory of species for structures associated to vector spaces over finite fields. 

This is not the first appearance of $q$-analogs nor the first use of vector spaces in the theory of species. D\'ecoste \cite{Decoste93,DecosteLabelle96} defined canonical $q$-counting series by means of $q$-substitutions in the cycle index series and in the asymmetry index series introduced by Labelle \cite{Labelle92}. Also, Joyal \cite{Joyal86} introduced the concept of a ``tensorial species,'' which is a functor from the category of finite sets (with bijections) to the category of vector spaces over a field of characteristic zero. However, the approach in this paper is different from the earlier work. First of all, we are not $q$-counting ordinary combinatorial structures but counting structures associated to vector spaces over the field of order $q$. Second, a \qsp is a functor from a category of vector spaces to the category of sets, whereas a tensorial species is a functor in the opposite direction.

First we recall the definition of a combinatorial species. Let $\mathbb{B}$ be the category whose objects are finite sets and whose morphisms are bijections. A species is a  functor $F: \B \ra \B$ \cite{Joyal81, BLL94}. For a finite set $U$, the set $F[U]$ is a collection of structures on the set $U$. 
Let $\fq$ be the finite field of order $q$. Define $\vq$ to be the category whose objects are finite dimensional vector spaces over $\fq$ and whose morphisms are the linear isomorphisms. 

\begin{definition}A \textbf{\qsp} (or \textbf{species of structures over} $\fq$) is a functor $F: \vq \ra \B$.
\end{definition}

Let $\fq^{(\N)}$ be the vector space of countable dimension whose elements are vectors $(a_1,a_2,\ldots)$ with a finite number of non-zero components. Let $e_1,e_2,\ldots$ be the standard basis and let $E_n$ be the span of $e_1,\ldots,e_n$. Then 
\[ E_0 \subset E_1 \subset \cdots  \subset E_n \subset E_{n+1} \subset \cdots \]
is an increasing sequence of subspaces whose union is $\fq^{(\N)}$. Let $\g_n$ be the order of the general linear group $\glnq$ of invertible $n \times n$ matrices over $\fq$ and define $\g_0=1$. Recall that
\[  \g_n = \prod_{i=0}^n (q^n-q^i)  .\]

\begin{definition} The \textbf{generating series} of a \qsp  $F$ is the power series
\[  \widehat{F}(x)= \sum_{n \geq 0} \frac{f_n}{\g_n} x^n, \]
where $f_n= | F[E_n] | $ is the number of elements in $F[E_n]$.
\end{definition}

\begin{example}   
The \qsp of \textbf{elements} $\eps$ is defined by $\eps[V]=V$ and $\eps[\phi]=\phi\,$ for an isomorphism $\phi : V \ra W$.
\[ \widehat{\eps}(x)=\sum_{n \geq 0}\frac{q^n}{\g_n} x^n .\]
\end{example}
\begin{example}  
The \qsp of \textbf{projective spaces} $\mathbb{P}$ with $\mathbb{P}[V]$ defined to be the set of one-dimensional subspaces of $V$. For an isomorphism $\phi:V \ra W$ and $L \in \mathbb{P}[V]$, we define  $\mathbb{P}[\phi](L) = \phi(L)$. The generating series is
\[ \widehat{\mathbb{P}}(x)= \sum_{n \geq 0} \frac{ [n]_q}{\g_n} x^n , \]
where $[n]_q=1 +q +q^2 + \cdots +q^{n-1}$ is the $q$-analog of $n$. (Note that $ [0]_q = 0$.)
\end{example}
\begin{example}  
The \qsp of \textbf{endomorphisms} is defined by $F[V]=\End_{\fq}(V)$ and $F[\phi](\alpha)= \phi \circ \alpha \circ \phi^{-1}$ for $\phi:V \ra W$ and $\alpha \in \End_{\fq}(V)$. The generating series is 
\[ \widehat{F}(x)= \sum_{n \geq 0} \frac{q^{n^2}}{\g_n} x^n . \]
\end{example}
\begin{example}   
The \qsp of \textbf{automorphisms} is defined by $F[V]=\Aut_{\fq}(V)$ and $F[\phi](\alpha)= \phi \circ \alpha \circ \phi^{-1}$ for $\phi:V \ra W$ and $\alpha \in \Aut_{\fq}(V)$. The generating series is 
\[ \widehat{F}(x)= \sum_{n \geq 0} \frac{\g_n}{\g_n} x^n =\frac{1}{1-x}. \]
\end{example}

\begin{example}   
For the \qsp of \textbf{ordered bases} we define $F[V]$ to be the set of $n$-tuples $(v_1,v_2,\ldots,v_n) \in V^n$ where $n=\dim V$ and the $v_i$ are a basis of $V$. For an isomorphism $\phi:V \ra W$ we define 
\[  F[\phi](v_1,\ldots,v_n)= (\phi(v_1),\ldots,\phi(v_n)) .\] The generating series is also  
\[  \widehat{F}(x)=\frac{1}{1-x}. \]
\end{example}

\begin{example}  
The \qsp of \textbf{vector spaces} $\mathcal{V}$ is defined by $\mathcal{V}[V]= \{V\}$ with generating series
\[ \widehat{\mathcal{V}}(x)=\sum_{n \geq 0} \frac{x^n}{\g_n}  .\]
The \qsp of \textbf{non-zero vector spaces} $\mathcal{V}_+$ is defined to be
\[ \mathcal{V}_+[V]= \begin{cases}
				\emptyset,  & \text{if $ \dim V=0$,}\\
				\{V\},  & \text{if $ \dim V > 0$.}
		 	    \end{cases}
\]	
with generating series
\[  	\widehat{\mathcal{V_+}}(x)=\sum_{n \geq 1} \frac{x^n}{\g_n}  .\]	    
\end{example}

\begin{example}  
The \qsp $F$ defined by $F[V]$ being the set of $k$-\textbf{dimensional subspaces of} $V$ has generating series \[ \widehat{F}(x)= \sum_{n \geq 0} \qbinom{n}{k}\frac{x^n}{\g_n} , \]
where 
\[  \qbinom{n}{k} = \frac{\g_n}{\g_k \g_{n-k} q^{k(n-k)}}  \]
is the $q$-binomial coefficient.  
\end{example}

\begin{example}  
Let $\Gamma$ be a group. The \qsp of \textbf{representations} is defined by $F[V]=\Hom(\Gamma,\Aut_{\fq}(V))$. For an isomorphism $\phi:V \ra W$ and a representation $\rho: \Gamma \ra \Aut_{\fq}(V)$. we have 
\[ F[\phi](\rho)= \phi \circ \rho \circ \rho^{-1}  .\]
The generating series depends on  $\Gamma$. Some results for finite groups are given in \cite{CTY00}. For cyclic groups one may also consult \cite{Morrison05seqmat}.
\end{example}

\begin{remark} Additional examples of generating series of \qsp are given in  \cite{Morrison05seqmat}. They include direct sum decompositions (splittings), flags of subspaces, linear and projective derangements, and diagonalizable, cyclic, or separable endomorphisms. 
\end{remark}

\begin{definition}
Two structures $s, t \in F[V]$ are \textbf{isomorphic}, indicated $s \sim t$, if there exists $\alpha \in \Aut_{\fq}(V)$ such that $F[\alpha](s)=t$. The number of isomorphism classes in $F[E_n]$ is denoted by $\widetilde{f}_n$ and the \textbf{type generating series} of $F$ is the formal power series
\[ \widetilde{F}(x)=\sum_{n \geq 0} \widetilde{f}_n x^n.  \]
\end{definition}

\begin{example}
The \qsp of ordered bases has only one isomorphism type in each dimension. Thus, $\widetilde{F}(x)= 1/(1-x)$. The \qsp of automorphisms, which  has the same generating series as the \qsp of ordered bases, has for $\widetilde{f}_n$ the number of conjugacy classes of invertible $n \times n$ matrices. It is shown in \cite{Morrison05seqmat} that the type generating series is
\[ \widetilde{F}(x) = \prod_{r \geq 1} \frac{1-x^r}{1-qx^r} . \]
\end{example}

In order to define the cycle index series of a \qsp we summarize the rational canonical form of a linear endomorphism. For $\sigma \in \End_{\fq}(V)$, $V$ is a module over $\fq[z]$ with $f(z)\cdot v$ defined to be $f(\sigma)(v)$. Then $V$ decomposes uniquely as  a direct sum of primary cyclic modules, which are modules of the form $\fq[z]/(\phi ^i)$ for some monic, irreducible polynomial $\phi$ and some positive integer $i$. Let $e_{\phi, i}(\sigma)$ be the number of copies of $\fq[z]/(\phi ^i)$ that occur in the decomposition of $V$. These integers are the invariants of $\sigma$ that completely describe its conjugacy class within $\End_{\fq}(V)$. There is a basis of $V$ for which the matrix representation of $\sigma$ is a block diagonal form consisting of $e_{\phi,i}$ copies of the companion matrix of $\phi^i$.
An endomorphism is an automorphism precisely when the polynomial $z$ does not occur among the invariants.
\begin{definition}  
The \textbf{cycle index series} of a \qsp $F$ is a formal power series in an infinite number of variables $x_{\phi,i}$ where $\phi$ ranges over the irreducible monic polynomials in $\fq[z]$, other than $\phi(z)=z$, and $i$ is a positive integer. We define this series to be
\[ \z_F=\sum_{n \geq 0} \frac{1}{\g_n} \sum_{\sigma \in \Aut(E_n)} \text{fix } F[\sigma]
               \prod_{\phi, i} x_{\phi,i}^{e_{\phi,i}(\sigma)},
\] 
where $\fix F[\sigma]$ is the number of fixed points of $F[\sigma]$.          
\end{definition}
\begin{remark}
The cycle index series defined here is not the same as those defined by Kung \cite{Kung81} and Stong \cite{Stong88}, although it bears a strong resemblance to them. 
\end{remark}

The cycle index series $\z_F$ can be specialized to give both the generating series $\wh{F}(x)$ and the type generating series $\wt{F}(x)$. In order to do so it is helpful to order the monic irreducible polynomials by putting $z-1$ first, then the rest of those of degree one, and then in order of increasing degree. Then we use the following equivalent notations:
\[ \z_F =\z_F((x_{\phi,i}))= \z_F((x_{z-1,1},x_{z-1,2},\ldots,x_{z-1,i},\ldots ),\ldots, (x_{\phi,1},x_{\phi,2},\ldots,x_{\phi,i},\ldots ), \ldots ).
\]
\begin{proposition} The generating series of $F$ is obtained from $\z_F$ by setting $x_{z-1,1}=x$ and by setting $x_{z-1,i}=0$ for $ i \geq 2$ and $x_{\phi,i}=0$ for all other $\phi$ and $i$. That is,
\[ \wh{F}(x)= \z_F((x,0,0,\ldots),(0,0,\ldots),\ldots,(0,0,\ldots),\ldots) . \]
\end{proposition}
\proof 
  \begin{align*} 
   \z_F((x,0,0,\ldots),(0,0,\ldots),\ldots,(0,0,\ldots),\ldots) 
            &= \sum_{n} \frac{1}{\g_n} \text{fix }F[I_n] \, x^n \\
            &= \sum_{n}\frac{f_n}{\g_n} \, x^n \\
            &= \wh{F}(x).
  \end{align*}
\qed

\begin{proposition} The type generating series $\wt{F}(x)$ is obtained from $\z_F$ by setting $x_{\phi,i}=x^i$ for all $\phi$ and $i$.
\end{proposition}
\proof
  \begin{align*}
  \z_F((x,x^2,\ldots,x^i,\ldots),(x,x^2,\ldots,x^i,\ldots),\ldots) 
            &= \sum_{n} \frac{1}{\g_n}  \sum_{\sigma \in \Aut(E_n)}\fix F[\sigma]\prod_{\phi,i}x^{i e_{\phi,i}(\sigma)}\\
            &= \sum_{n}\frac{1}{\g_n}\sum_{\sigma \in \Aut(E_n)}\fix F[\sigma] x^{\sum_{\phi,i}i e_{\phi,i}(\sigma)}\\
            &=\sum_{n}\frac{1}{\g_n}\sum_{\sigma \in \Aut(E_n)}\fix F[\sigma] \, x^n \\
            &= \sum_{n} \wt{f}_n \,x^n .
\end{align*}
The last step uses  Burnside's Lemma for the number of orbits of a finite group action. In this case the orbits of $\Aut(E_n)$ acting on $F[E_n]$ are the isomorphism classes of structures in $F[E_n]$.
\qed

\begin{definition}
Two \qsp $F$ and $G$ are \textbf{isomorphic} \qsp if they are isomorphic as functors, i.e. there exists an invertible morphism of functors (natural transformation) $\eta: F \ra G$. We consider isomorphic \qsp to be equal and write $F = G$. (The term ``combinatorially equal'' is used by Bergeron, Labelle and Leroux \cite{BLL94}.)
\end{definition}

\section{Sums and Products}  
\begin{definition}
Given \qsp $F$ and $G$ we define their \textbf{sum}  \qsp $F+G$ by
\begin{align*}  (F+G)[V]  &= F[V] \uplus G[V] \quad \textrm{(disjoint union)} \\
			(F + G)[\phi](s) &= 
			     \begin{cases}
			          F[\phi](s),  & \text{if $ s \in F[V] $,}\\
			          G[\phi](s),  & \text{if $ s \in G[V]$.}
			        \end{cases}  
\end{align*}
The \textbf{product} \qsp $F \cdot G$ is defined on objects by			        
\[  (F \cdot G)[V] =  \bigcup_{V_1 \oplus V_2 = V} F[V_1] \times G[V_2]  .\]
For an isomorphism $\phi: V \ra W$ and for  $(s,t) \in F[V_1] \times G[V_2]$ we define 
\[ (F\cdot G)[\phi](s,t)=(F[\phi_1](s),G[\phi_2](t)) , \]
where $\phi_i = \phi | V_i$. 
\end{definition}

\begin{proposition} For \qsp $F$ and $G$ the generating, type generating, and cycle index series of their sum and product satisfy
\begin{align*}
  \widehat{(F+G)}(x) &= \widehat{F}(x) + \widehat{G}(x)  \\
  \widetilde{(F+G)}(x) &= \widetilde{F}(x) + \widetilde{G}(x) \\
  \z_{F+G} &= \z_F + \z_G \\  \\
  \widehat{(F \cdot G)}(x) &= \widehat{F}(x) \cdot \widehat{G}(x)  \\
  \widetilde{(F \cdot G)}(x) &= \widetilde{F}(x) \cdot \widetilde{G}(x) \\
   \z_{F \cdot G} &= \z_F \cdot \z_G.
\end{align*}
\end{proposition}
\proof We prove only the assertions about the product, since those for the sum are straightforward. Let $H = F\cdot G$. Then $h_n$,   the cardinality of $H[E_n]$, is given by 
\[  h_n = \sum_{0 \leq k \leq n}  \sum_{\substack{ \dim V_1 = k \\ \dim V_2=n-k \\ V_1 \oplus V_2=E_n} }f_k g_{n-k}  .\]
There are $\g_n / \g_k \g_{n-k} $ direct sum decompositions $E_n=V_1 \oplus V_2$ where $\dim V_1 =k$, and so
\[ h_n= \sum_{k=0}^{n} \frac{\g_n}{\g_k \g_{n-k}}f_k g_{n-k} .\] 
Therefore, 
\begin{align*}
 \sum_{n \geq 0} \frac{h_n}{\g_n}x^n  &=
     \sum_{n \geq 0} \sum_{k=0}^{n} \frac{f_k}{\g_k }\frac{g_{n-k}}{\g_{n-k}}x^n \\
      \widehat{H}(x) &= \widehat{F}(x) \widehat{G}(x) .
\end{align*}

To prove that $\wt{H} = \wt{F}\wt{G}$ we need to prove that 
\[ \widetilde{h}_n = \sum_{k=0}^n \widetilde{f}_k \widetilde{g}_{n-k} .\]
This will follow by showing that there is a bijection
\be \bigcup_{k=0}^n  F[E_k]/ \!\! \sim  \times  \, \,G[E_{n-k}]/\! \! \sim \, \longrightarrow \, H[E_n]/\! \! \sim.
\label{bijection}
 \ee
Let $\iota$ be the isomorphism from $E_{n-k}$ to the subspace of $E_n$ spanned by  $e_{k+1}, \ldots, e_n$ defined by $e_i \mapsto e_{k+i}$. Then in (\ref{bijection}) we  map the pair of isomorphism classes $([s],[t])$ to the isomorphism class of $(s,G[\iota](t))$ in $H[E_n]$. It is routine to see that the map is well-defined and injective. To see that it is surjective consider a structure in $H[E_n]$, say $(s,t)$ where $s \in F[V_1]$,  $t \in G[V_2]$, and $V_1 \oplus V_2 = E_n$ with $\dim V_1=k$. Choose $\alpha \in \Aut(E_n)$ such that $\alpha$ maps $V_1$ to $E_k$ and $V_2$  to $\iota(E_{n-k})$. Then $(s,t)$ is isomorphic to
$H[\alpha](s,t) \in F[E_k] \times G[\iota(E_{n-k})]$ and so every isomorphism class in $H[E_n]$  is the image of a pair from $F[E_k] \times G[E_{n-k}]$. The map defined by (\ref{bijection}) is a bijection.

For the final claim that we begin with the definition
\[  \z_{F \cdot G} = \sum_{n \geq 0} \frac{1}{\g_n} \sum_{\sigma \in \Aut(E_n)} \text{fix }( F\cdot G)[\sigma]
               \prod_{\phi, i} x_{\phi,i}^{e_{\phi,i}(\sigma)} .\]
A structure $(s,t)$, with $s \in F[V_1]$, $t \in G[V_2]$,$ V_1 \oplus V_2 = E_n$, is fixed by $F \cdot G$ if and only if $\sigma_i = \sigma | V_i$ is an automorphism of $V_i$ for $i=1,2$ and $s$ is fixed by $F[\sigma_1]$ and $t$ is fixed by $G[\sigma_2]$. Thus, $e_{\phi,i}(\sigma)=e_{\phi,i}(\sigma_1)+e_{\phi,i}(\sigma_2)$ and
\[  \fix (F \cdot G)[\sigma]= \sum_{V_1 \oplus V_2=E_n} \sum_{\sigma_1 \in \Aut(V_1)}
     \sum_{\sigma_2 \in \Aut(V_2)} \fix F[\sigma_1] \, \fix G[\sigma_2] .\]      
We group the terms on the right according to $m=\dim V_1$.  Recall that there are $\g_n/\g_m \g_{n-m}$ decompositions of $E_n$ into a direct sum of subspaces of dimension $m$ and $n-m$. This gives
\[  \fix (F \cdot G)[\sigma]= \sum_{m=0}^n  \frac{\g_n}{\g_m \g_{n-m}}\sum_{\sigma_1 \in \Aut(E_m)}
         \sum_{\sigma_2 \in \Aut(E_{n-m})}  \fix F[\sigma_1] \fix[G\sigma_2]  .\]  
Therefore,
\begin{align*}  \z_{F \cdot G} &= \sum_{n \geq 0} 
                                  \sum_{m=0}^n  \frac{1}{\g_m \g_{n-m}}\sum_{\sigma_1 \in \Aut(E_m)}
                                  \sum_{\sigma_2 \in \Aut(E_{n-m})}  \fix F[\sigma_1] \fix[G\sigma_2] 
                                  \prod_{\phi, i} x_{\phi,i}^{e_{\phi,i}(\sigma_1)+e_{\phi,i}(\sigma_2)}  \\
                               &=\left(  \sum_{m \geq 0} \sum_{\sigma_1 \in \Aut(E_m)} \fix F[\sigma_1]
                                   \prod_{\phi,i}x_{\phi,i}^{e_{\phi,i}(\sigma_1)} \right)
                                   \left(  \sum_{k \geq 0} \sum_{\sigma_2 \in \Aut(E_k)} \fix G[\sigma_2]
                                   \prod_{\phi,i}x_{\phi,i}^{e_{\phi,i}(\sigma_2)} \right)  \\
                                &=\z_F \cdot \z_G
\end{align*}
\qed

\begin{remark} The isomorphism classes of \qsp form a commutative semi-ring using the sum and product operations. The additive and multiplicative identities $0$ and $1$ are defined by
\begin{align*}
  0[V] &= \emptyset  \\
  1[V] &= \begin{cases}
                \{V\}, & \text{if $\dim V=0$,} \\
                \emptyset, & \text{if $\dim V > 0$.}
             \end{cases}
\end{align*}
The \qsp $n$, defined to be the $n$-fold sum $1 + 1 + \cdots + 1$, is the \qsp that has exactly $n$ structures on $V=\{0\}$ and none on any vector space $V$ of positive dimension.  Thus, the natural numbers are embedded in the semi-ring of isomorphism classes of \qsp.
The associated ring constructed from formal differences is the \textbf{ring of virtual \qsp}. 
\end{remark}

\section{Symmetric Powers and Assemblies}
For a \qsp $F$ we let $F^n$ denote the $n$-fold product of $F$ with itself. There is a natural action of the symmetric group $S_n$ on $F^n$ permuting the components of a structure $(s_1,\ldots,s_n) \in F^n[V]$.
\begin{definition}
Let $F$ be a \qsp and $n$ a positive integer.  Define the \qsp  $F^\sym{n}$, the \textbf{$n$th-symmetric power} of $F$, by \[ F^\sym{n}[V] = F^n[V]/S_n .\]
 A structure in $F^\sym{n}[V]$ is a multi-set  $\{ s_1,s_2,\ldots,s_n \}$.
 \end{definition} 
\begin{proposition}
If $F[0]=\emptyset$, then $\wh{F^\sym{n}}(x)=\wh{F}(x)^n/n! $.
\end{proposition}
\proof  With the hypothesis that $F[0]= \emptyset $, any structure $\{ s_1,\ldots, s_n \}$ has an associated direct sum decomposition $V_1 \oplus \cdots \oplus V_n$ in which none of the $V_i$ is the zero subspace. Therefore, all the $V_i$ are distinct subspaces and the action of $S_n$ on $F^n[V]$ is free. It follows that the  cardinality of $F^\sym{n}[V]$ is the cardinality of $F^n[V]$ divided by $n!$. This means that each coefficient of the generating series for $F^\sym{n}$ is obtained from the corresponding coefficient of the generating series for $F^n$ by dividing by $n!$. 
\qed

\begin{remark}
The presence of zero subspaces in direct sum decompositions of $V$ complicates the counting of the structures in $F^\sym{n}[V]$. In order to construct symmetric powers without allowing trivial subspaces in the decompositions, one may use the symmetric powers of the \qsp $F_{+}$, which is the same as $F$ in positive dimensions but has no structures on the zero vector space.
\end{remark}

\begin{definition}
Let $F$ be a \qsp. We call a structure $\{ s_1, s_2,\ldots,s_n\} \in F^\sym{n}[V]$ an \textbf{assembly of $F$-structures on} $V$ if $s_i \in F[V_i]$ and $V=V_1 \oplus \cdots \oplus V_n$ is a non-trivial direct sum decomposition of $V$. (Non-trivial means that none of the subspaces is zero.) \end{definition}

\begin{theorem} \label{symsum}
For \qsp $F$ and $G$ there is an isomorphism of  \qsp
\[  (F +G)^{\sym{n}} = \sum_{m=0}^{n}F^{\sym{m}} \cdot G^\sym{n-m} .\]
\end{theorem}
\proof  
An assembly in $(F +G)^{\sym{n}} [V]$ is a set of structures $\{s_1,\ldots,s_m,t_1,\ldots,t_{n-m} \}$ where $s_i \in F[U_i]$ and $t_i \in G[W_i]$ for a splitting $V=U_1 \oplus \cdots \oplus U_m \oplus W_1 \oplus \cdots \oplus W_{n-m} $. Such an assembly is a structure in $(F^{\sym{m}} \cdot G^{\sym{n-m}})[V]$ with the decomposition $V=V_{1} \oplus V_{2}$ where $V_{1}=U_{1}\oplus \cdots \oplus U_m$ and $V_{2}= W_1 \oplus \cdots \oplus W_{n-m}$.
\qed
\begin{definition}
Let $F$ be a \qsp with $F[0]=\emptyset$. Define  $E \circ F= 1+\sum_{n \geq 1}F^\sym{n}$, the \qsp of \textbf{assemblies of $F$-structures}. We also use the notation $ E(F) = E \circ F$.
\end{definition}

\begin{remark}
In the setting of combinatorial species there is a general notion of ``substitution'' or ``partitional composition'' for two species, and assemblies of structures are a special case involving the species of sets $E$.  See \cite{BLL94} for more information and an explanation of the notation. It is possible to define the substitution $H \circ F$ for a combinatorial species $H$ (a functor from the category $\mathbb{B}$ to itself) and a \qsp $F$. Letting $H=E_{n}$ be the species of $n$-sets, we have $E_{n}\circ F=F^\sym{n}$ when $F[0]=\emptyset$, but we have no more interesting examples for $H$ is anything other than $E$ and $E_{n}$, and so we do not follow that thread here.
Symmetric powers such as defined here and further generalizations are considered by Joyal \cite{Joyal86} in the theory of ordinary species.
\end{remark}

\begin{theorem} For \qsp $F$ and $G$ there is an isomorphism of \qsp
\[     E \circ (F + G) = (E \circ F) \cdot  (E \circ G).  \]
\end{theorem}
\proof 
{}From Theorem \ref{symsum} we see that
\begin{align*}  E \circ (F + G)  &=  1 + \sum_{n \geq 1}(F + G)^{\sym{n}}  \\
                               &= 1 + \sum_{n \geq 1}\sum_{m=0}^{n} F^{\sym{m}} \cdot G^\sym{n-m}  \\
                               &= \left(1 + \sum_{m \geq 1}F^{\sym{m}} \right) \cdot
                                    \left( 1 + \sum_{k \geq 1}G^{\sym{k}} \right) \\
                               &= (E \circ F) \cdot (E \circ G)     
\end{align*}
\qed

\begin{theorem}
The generating series for $E \circ F$ is given by
\[ (\wh{E \circ F})(x) = \exp( \wh{F} (x) )  .\]
\end{theorem}
\proof 
{}From the definition $E \circ F= \sum_{n \geq 1}F^\sym{n}$ we have
\begin{align*}   
	(\wh{E \circ F})(x) &= 1 + \sum_{n \geq 1}\wh{F^{\sym{n}}}(x) \\
				&=  1  + \sum_{n \geq 1}\frac{1}{n!}\wh{F}(x)^{n} \\
				&= \exp( \wh{F} (x))  .
\end{align*}
 \qed

\begin{theorem}
 Let $F$ be a \qsp and suppose $F[0]=\emptyset$. Then the type generating series of $E \circ F$ is given by
 \[  (\wt{E \circ F})(x) = \prod_{n \geq 1}\frac{1}{(1-x^{n})^{\wt{f}_{n}}}  \]
 or, alternatively, by
 \[   (\wt{E \circ F})(x) =\exp \left( \sum_{n \geq 1}\frac{1}{n}\wt{F}(x^{n})\right)  .\]
 \end{theorem}
 
 \proof
 Decompose $F$ as the sum $\sum_{n \geq 1}F_{n}$, where 
 \[  F_{n}[V]= \begin{cases}  F[V],   & \text{if $ \dim V=n$,}\\ 
                                           \emptyset, &\text{otherwise.}
                      \end{cases}
\] 
Then 
\begin{align*}
   \wt{ E \circ F} &= \wt{E \circ \sum_{n \ge 1}F_{n} } \\
                 &= \prod_{n \geq 1}( \wt{E \circ F_{n}}) .
\end{align*}
In order to find $( \wt{E \circ F_{n}})(x)$ we observe that the isomorphism types of assemblies of $F_{n}$-structures on a vector space $V$ (whose dimension must be $kn$ for some $k \geq 1$) correspond with multsets of size $k$ chosen from a set of size $\wt{f}_{n}$. There are $\binom{\wt{f}_{n}+k-1}{k}$ such multisets. Therefore,  $( \wt{E \circ F_{n}})(x)$ is the generating function
\[  \frac{1}{(1-x^{n})^{\wt{f}_{n}}} = \sum_{k \geq 0}\binom{\wt{f}_{n}+k-1}{k} x^{k} .\]
For the second formula of the theorem, take the $\log$ of the product, expand the $\log$ of each term as a series, and switch the order of summation.
\qed

\begin{remark}  In the terminology of Wilf \cite{Wilf94} the first formula gives the hand enumerator for a prefab. The second formula shows exactly the same form as for $(\wt{E \circ F})$ in the case of combinatorial species \cite[p. 43]{BLL94}. In that case
\[  (\wt{E \circ F})(x)= Z_{E}(\wt{F}(x), \wt{F}(x^{2}),\wt{F}(x^{3}),\ldots) \]
and
\[  Z_{E}(x_{1},x_{2},x_{3},\ldots)= \exp \left( x_{1} + \frac{x_{2}}{2} + \frac{x_{3}}{3},\ldots  \right). \]
 \end{remark}

\begin{example}  \label{unosplit}
 The \qsp of unordered splittings (i.e. direct sum decompositions in which the order of the summands is not important) is given by $E \circ \mathcal{V}_{+}$. A structure in $( E \circ \mathcal{V}_{+})[V]$ is $\{ V_{1},\ldots,V_{n} \}$ where $V= V_{1} \oplus \cdots \oplus V_{n}$ and $\dim V_{i} > 0$. The generating series is
 \[ ( \wh{E \circ \mathcal{V}_{+}})(x)=\exp \left( \sum_{n \geq 1}\frac{x^{n}}{\g_{n}} \right). \]
The coefficients of this series, which count the number of splittings, give $q$-analogs of the Bell numbers, which count the number of set partitions \cite{Morrison04}.
The type generating series is
\[  (\wt{E \circ \mathcal{V}_+} )(x)= \prod_{n \geq 1} \frac{1}{1-x^n}  , \]
which is the partition generating function. For a vector space $V$ of dimension $n$ the isomorphism type of a splitting is completely determined by the dimensions of the direct summands, i.e. by a partition of $n$.
\end{example}

\begin{example}
 Let $D$ be the \qsp for which $D[V]$ is the set of diagonalizations of endomorphisms of $V$. Thus, a structure in $D[V]$ is a pair $(\alpha, \{ V_{1},V_{2},\ldots,V_{n} \} )$ where $\alpha \in \End(V)$, $V=V_{1}\oplus \cdots \oplus V_{n}$, $\dim V_{i}=1$, and $\alpha | V_{i} = \lambda_{i}I$ for scalars $\lambda_{i} \in \fq$. Then $D=E \circ \mathbb{F}$ where $\mathbb{F}$ is the \qsp defined
\[  \mathbb{F}[V] = 
	\begin{cases}  
		\fq, &\text{if  $ \dim V=1$, }\\
		\emptyset,  &\text{otherwise.}
	\end{cases}
\]
Then
\begin{align*}  
	\wh{\mathbb{F}}(x) &=\frac{q}{q-1}x \\
	\wh{D}(x) &= \exp\left( \frac{q}{q-1}x \right)  \\
	\wt{D}(x) &= \frac{1}{(1-x)^{q}} .
\end{align*}
\end{example}

\begin{example}
Let $\fq^{\times}$ be the multiplicative group of $\fq$. Modify the previous example by defining 
\[\mathbb{F}^{\times}[V]= 
	\begin{cases}
		\fq^{\times}, &\text{if  $ \dim V=1$, }\\
		\emptyset,  &\text{otherwise.}
	\end{cases}
\]
Then an assembly of $\mathbb{F}^{\times}$-structures corresponds to a diagonalization of an automorphism. Let $D^{\times}= E \circ \mathbb{F}^{\times}$. We have the generating series
\begin{align*}
	\wh{\mathbb{F}^{\times}}(x) & = x \\
	\wh{D^{\times}}(x) &= e^{x}  \\
	\wt{D^{\times}}(x) &= \frac{1}{(1-x)^{q-1}}.
\end{align*}
\end{example}

\begin{remark}
In \cite{Morrison04} assemblies of structures of \qsp are treated from the point of view of ``$q$-exponential families,'' which are analogs of the exponential families in Wilf's book \emph{generatingfunctionology} \cite{Wilf94}. The formula $(\wh{E \circ F})(x)=\exp(\wh{F}(x))$ is the exponential formula giving the one-variable hand enumerator. The formula for $(\wt{E \circ F})(x)$ is the one-variable hand enumerator for a ``prefab'' \cite[Theorem 3.14.1]{Wilf94}. The two-variable hand enumerators for exponential families and prefabs require the use of \emph{weighted} \qsp for their statments.
\end{remark}
\begin{question} 
Is there a formula for $\mathcal{Z}_{E \circ F}$ in terms of $\mathcal{Z}_{F}$ as there is for ordinary combinatorial species? It is not clear what to expect even for the most basic example of a \qsp $F$ with exactly one structure in dimension one.
\end{question}
\begin{question}
Is there a formula for $\mathcal{Z}_{F^{\sym{n}}}$ in terms of $\mathcal{Z}_{F}$? 
\end{question}
\begin{question}
Is there a formula for $\wt{ F^{\sym{n}}}(x)$  in terms of $\wt{ F}(x)$?  An answer to the previous question should give an answer to this one, but it may be more efficient to bypass the cycle index.
 In principle, the formula for $( \wt{E \circ F})(x)$ could be obtained from 
\[  ( \wt{E \circ F})(x) =  1+ \sum_{n \geq 1} \wt{ F^{\sym{n}}}(x) . 
 \]
\end{question} 
 
\section{Weighted \textit{q}-Species}
\begin{definition} Let $\A$ be a commutative ring. An $\A$\textbf{-weighted set} is a pair $(A,w)$ where $A$ is a set and $w:A \ra \A$ is a function called a \textbf{weighting}.
\end{definition}
We let $\B_{\A}$ be the category whose objects are finite $\A$-weighted sets and whose morphisms are bijections that commute with the weightings.
\begin{definition}
An $\A$\textbf{-weighted \qsp} is a functor $F_{w}: \vq \ra \B_{\A}$.
\end{definition}

\begin{example} The \qsp of ordered or unordered splittings can be made into a weighted \qsp by defining the weight of a splitting to be $t^{k} \in \Q[t]$ where $k$ is the number of summands in the decomposition.
\end{example}

In the category of weighted sets, the cardinality of a set is replaced by the \textbf{inventory} or \textbf{total weight}
\[  |A|_{w} = \sum_{a \in A}w(a) .\]
 For weighted sets $(A,w)$ and $(B,v)$ their sum is defined to be $(A+B, \mu)$ where $A+B$ is the disjoint union of $A$ and $B$ and 
 \[ \mu(x)= \begin{cases}   w(x), &\text{if  $x \in A$, }\\
				  v(x),  & \text{if $x \in B$.}
		\end{cases}
\]
Their product is defined to be $(A \times B, \rho)$	
where 
 \[ \rho(a,b)=w(a)v(b).\]	
 One easily checks that
 \begin{align*}
   |A+B|_{\mu} &= |A|_{w} + |B|_{v}	\\
   |A\times B |_{\rho} &= |A|_{w} \, |B|_{v}.
\end{align*} 
With the trivial weighting $w(a)=1$ any set is a $\Z$-weighted set, $\Z$ being the ring of integers.
\begin{definition}
For an $\A$-weighted \qsp $F=F_{w}$ the \textbf{generating series} of $F$ is the power series with coefficients in $\A$ given by
\[  \wh{F_{w}}(x)=\sum_{n \geq 0} |F[E_{n}]|_{w} \frac{x^{n}}{\g_{n}} .\]

The \textbf{type generating series} of $F$ is the power series
\[   \wt{F_{w}}(x)=\sum_{n \geq 0} |F[E_{n}]/\sim |_{w} x^{n} .\]
The set of isomorphism classes $F[E_{n}]/\sim$ is a weighted set with the weight of a class being the weight of any representative. They are all equal since weights are preserved by morphisms in the category of weighted sets.

The \textbf{cycle index series} of $F$ is defined by
\[ \z_{F_{w}} = \sum_{n \geq 0} \frac{1}{\g_n} \sum_{\sigma \in \Aut(E_n)} |\Fix F[\sigma]|_{w}
               \prod_{\phi, i} x_{\phi,i}^{e_{\phi,i}(\sigma)}.
\]
The fixed point set  $\Fix F[\sigma]$ inherits the weighting as a subset of $F[E_{n}]$.            
\end{definition}
\begin{proposition}
Let $F_{w}$ and $G_{v}$ be $\A$-weighted \qsp. Then
\begin{align*}
  \widehat{(F_w+G_u)}(x) &= \widehat{F_w}(x) + \widehat{G_u}(x)  \\
  \widetilde{(F_w+G_u)}(x) &= \widetilde{F_w}(x) + \widetilde{G_u}(x) \\
  \z_{F_w+G_u} &= \z_{F_w} + \z_{G_u} \\  \\
  \widehat{(F_w \cdot G_u)}(x) &= \widehat{F_w}(x) \cdot \widehat{G_u}(x)  \\
  \widetilde{(F_w\cdot G_u)}(x) &= \widetilde{F_w}(x) \cdot \widetilde{G_u}(x) \\
   \z_{F_w\cdot G_u} &= \z_{F_w }\cdot \z_{G_u}.
\end{align*}
\end{proposition}
\proof   The proof is similar to the proof for the unweighted case.  \qed

\begin{corollary} If $F_{w}$ is a weighted \qsp, then $F_{w}^{n}$ is weighted with 
\[ \wh{F_{w}^{n}}(x) = (\wh{F_{w}}(x))^{n} . \]
Furthermore, if $F_{w}[0]=\emptyset$, then
\[\wh{F_{w}^{\sym{n}}}(x) = (\wh{F_{w}}(x))^{n}/n! \]
and the generating series for the weighted \qsp of assemblies of $F_{w}$-structures is
\[ \wh{(E \circ F_{w}})(x)= \exp \left( \wh{F_{w}}(x) \right) .\]
\end{corollary}

\begin{example} 
The \qsp of unordered splittings is the \qsp of assemblies $E \circ \mathcal{V}_{+}$ from Example 4.3. The weighting defined to be $t^{k}$ on a splitting with $k$ direct summands is actually the induced weighting coming from  the weighting for $\mathcal{V}_{+}$ given by $w(\{V\})=t \in \Q[[t]]$. The generating series for the weighted $\mathcal{V}_{+}$ is 
\[ \wh{\mathcal{V}_{+}}(x)= \sum_{n \geq 1} t \frac{x^{n}}{\g_{n}} ,\] and so the generating series for the weighted \qsp of splittings is 
\[ \exp \left( t \sum_{n \geq 1} \frac{x^{n}}{\g_{n}}\right) .\]
( In the language of \cite{Morrison04} this is the ``exponential formula'' for the  two variable hand enumerator for the $q$-exponential family of splittings.)
\end{example}

\providecommand{\bysame}{\leavevmode\hbox to3em{\hrulefill}\thinspace}
\providecommand{\MR}{\relax\ifhmode\unskip\space\fi MR }
\providecommand{\MRhref}[2]{%
  \href{http://www.ams.org/mathscinet-getitem?mr=#1}{#2}
}
\providecommand{\href}[2]{#2}

 \end{document}